\documentclass[12pt]{article}
\usepackage{amsmath}
\usepackage{amssymb}
\usepackage{epsf}
\usepackage[all]{xy}
\usepackage[english,francais]{babel}

\input xy
\xyoption{all}

\newtheorem{theor}{Theorem}[section]
\newtheorem{propo}[theor]{Proposition}
\newtheorem{coro}[theor]{Corollary}

\newtheorem{defpropo}[theor]{Definition-Proposition}
\newtheorem{rmke}[theor]{Remark}
\newtheorem{lemmae}[theor]{Lemma}

\newenvironment{theo}{\begin{theor} \sf}{\end{theor}}
\newenvironment{prop}{\begin{propo} \sf}{\end{propo}}
\newenvironment{cor}{\begin{coro} \sf}{\end{coro}}
\newenvironment{defprop}{\begin{defpropo} \sf}{\end{defpropo}}
\newenvironment{rmk}{\begin{rmke} \normalfont}{\end{rmke}}
\newenvironment{lemma}{\begin{lemmae} \sf}{\end{lemmae}}

\newenvironment{proof}{{\flushleft \bf Proof: \,} }

\newenvironment{ex}{{\flushleft \bf Example: \, }}

\def\data{{\flushleft \bf Data: \, }}
\def\hcm{HC_*^{\pi, \alpha, \beta} }
\def\hcme{HC_*^{\pi, \varepsilon, \varepsilon} }

\def\kk{ _{\beta}  k_{\alpha}}
\def\sp{S_{\pi}}
\def\ot{\otimes}

\def\lan{\Lambda _n}
\def\m{\frak{m}}

\advance\oddsidemargin by- 2.1cm
\advance\topmargin by- 1.2cm 
\numberwithin{equation}{section}

\title{\textsc{Cyclic homology of Hopf algebras}}

\author{Rachel Taillefer\thanks{Laboratoire G.T.A., CNRS UMR 5030, D\'epartement de Math\'ematiques CC 51, Universit\'e Montpellier II, 34095 Montpellier Cedex 5. 
email: taillefr\at math.univ-montp2.fr}}

\date{}

\begin{document}

\maketitle

%---------------------------------------------------------------------------%

%\input{absfcm.tex}

\selectlanguage{english}
\begin{abstract} A cyclic cohomology theory adapted to Hopf algebras has been introduced recently by Connes and Moscovici. In this paper, we consider this object in the homological framework, in the spirit of Loday-Quillen (\cite{LQ}) and Karoubi's work on the cyclic homology of associative algebras.  In the case of group algebras, we interpret the decomposition of the classical cyclic homology of a group algebra (\cite{B}, \cite{KV}, \cite{L}) in terms of this homology. We also compute both cyclic homologies for truncated quiver algebras.
\end{abstract}


\section{Introduction}

In this paper, we study a cyclic homology theory for  Hopf algebras, and we specialize our considerations mainly to group algebras and truncated quiver algebras.

Connes and Moscovici introduced cyclic cohomology of Hopf algebras
in~\cite{CM2}. They used it to solve a problem in noncommutative
geometry, namely the computation of the index of transversally
elliptic operators on foliations. They needed to compute an index pairing map associated to some convolution algebras $\mathcal{A}$ (see also~\cite{CM3}). To do this, they constructed some Hopf algebras $\mathcal{H} (n),$ endowed with  extra data, and they defined a cyclic cohomology for these, subject to the following condition: there should be a characteristic map from the cyclic cohomology of the $\mathcal{H} (n)$ to the classical cyclic cohomology of $\mathcal{A}$ for each trace $\tau: \mathcal{A} \rightarrow \mathbb{C}$ satisfying some invariance condition (see~\cite{CM2}, \cite{CM3}).

In \cite{Cr}, Crainic gives an alternative definition of this cyclic
cohomology in terms of $X-$complexes. He also computes some examples,
and constructs a non-commutative Weil complex which is related to the
cyclic cohomology of Hopf algebras. Cyclic homology of Hopf algebras
also appears in the final section of his thesis. It was also used by Gorokhovsky in~\cite{G} to construct some cyclic cocycles associated to a specific vector bundle.

In this paper, we consider  the homological framework, inspired by
Loday-Quillen and Karoubi's approach to the usual cyclic homology of
algebras obtained from the cohomological theory considered by
Connes. We proceed as follows: in the first section, we define a
cyclic homology of Hopf algebras, and describe Connes' long periodic
exact sequence (the SBI exact sequence) in this situation. Next, we study the case of group algebras, and we interpret the decomposition of the classical cyclic homology of group algebras, established by Burghelea (\cite{B}) and Karoubi and Villamayor (\cite{KV}), in terms of Connes and Moscovici's homology. We also compute the cyclic homology of the cyclic group algebras (as Hopf algebras). 

Finally, we consider the family of truncated quiver algebras. The Hochschild homology of a truncated quiver algebra with coefficients in itself was computed by Sköldberg in \cite{S}; this enables us to compute its classical cyclic homology. We then specialize to the Taft algebras, which are non-commutative, non-cocommutative Hopf algebras (of finite representation type), and compute their cyclic homology as Hopf algebras.

Unless specified, $k$ is a commutative ring.

{\bf Acknowlegements:} I am grateful to my advisor Claude Cibils and
to Jean-Michel Oudom for the many discussions we have had.

\section{Connes and Moscovici cyclic (co)homology }

In~\cite{CM}, A. Connes and H. Moscovici have defined a cyclic cohomology theory for Hopf algebras endowed with some specific data, including a grouplike element and a character satisfying some extra properties. This can be slightly extended to include a second grouplike element; we shall give here the dual version of this definition (cyclic homology of a Hopf algebra).

\data Let $k$ be a commutative ring, and let $H$ be a Hopf algebra
over $k.$ Let $\mu$ denote the multiplication, $\eta$ the unit,
$\Delta$ the comultiplication, $\varepsilon$ the counit, and $S$ the
antipode of $H.$ Suppose there are two characters $\alpha, \  \beta~: H \rightarrow k$ on $H,$ and a grouplike element $\pi$ in $H.$ Set $S_{\pi}= \pi S,$ and assume the following properties are satisfied: $$\begin{cases}
(i) \ \alpha (\pi) =1= \beta (\pi), \\
(ii) \ (\alpha * S_{\pi} * \beta)^2 = id, 
\end{cases} \ \ \ \label{cmcond}\textsc{(cm)}$$ where $*$ is the
convolution product (defined for linear maps from $H$ to $H$ or $k$).
The second identity can be rewritten with Sweedler's notation:
$$ \alpha (a^{(1)}) \alpha ( S ( a^{(4)})) \pi S^2 (a^{(3)}) \pi ^{-1} \beta ( S (a^{(2)})) \beta (a^{(5)})=a, \ \forall a \in H.$$

\begin{defprop} The following data define a cyclic module $C_*^{(\pi, \alpha, \beta)}(H):$ 
$$C_n^{(\pi, \alpha, \beta)}(H) = H^{\otimes n} \mbox{ for $n\geq 0,$}$$ and
$$\begin{array}{rcl}
d_i: H^{\otimes n+1}& \longrightarrow & H^{\otimes n}\\
a_0 \otimes \ldots \otimes a_n & \mapsto & \begin{cases} \alpha (a_0) a_1 \otimes \ldots \otimes a_n \mbox{ if $i=0,$} \\
a_0 \otimes \ldots \otimes a_{i-1} a_i \otimes \ldots \otimes a_n \mbox{ if $1 \leq i \leq n,$} \\
a_0 \otimes \ldots \otimes a_{n-1} \beta(a_n) \mbox{ if $i=n+1,$} \end{cases} \\
\\
s_i: H^{\otimes n}& \longrightarrow & H^{\otimes n+1}\\
a_1 \otimes \ldots \otimes a_n & \mapsto & \begin{cases} a_1  \otimes \ldots \otimes 1 \otimes a_i  \otimes \ldots \otimes a_{n} \mbox{ if $1 \leq i \leq n,$} \\
a_1  \otimes \ldots \otimes a_{n} \otimes 1 \mbox{ if $i=n+1,$ and} \end{cases} \\
\\
t_n: H^{\otimes n}& \longrightarrow & H^{\otimes n}\\
a_1 \otimes \ldots \otimes a_n & \mapsto & \alpha (a_1^{(1)} \ldots a_n^{(1)}) S_{\pi} (a_1^{(2)} \ldots a_n^{(2)}) \otimes a_1^{(3)} \otimes \ldots \otimes a_{n-1}^{(3)} \beta (a_n^{(3)}),
\end{array}$$ 
when $n \geq 1.$ When $n=0,$ the maps become: $d_0 = \alpha, \ d_1 = \beta, \  s_0= \eta$ (the unit of $H$), and $t_0 = id_k.$ 

The cyclic homology of $H$ is the homology of the usual bicomplexes associated to a cyclic module, or, when $k$ contains $\mathbb{Q}$, of the Hochschild complex factored by the cyclic action $t_.$ (see~\cite{L}). We shall denote this homology by $\hcm (H).$
\end{defprop}

\begin{proof}Most of the relations which should be satisfied by a
cyclic module are easy to check; the main difficulty arises for the
relation $t_n^{n+1} =id.$ We will use repeatedly the following property of $\sp :$ $\sp (ab)=\sp (b)S(a).$  Let us first compute the square of $t_n:$

\begin{eqnarray*} 
t_n^2 (a_1 \otimes \ldots \otimes a_n) &=& \alpha (a_1^{(1)} \ldots
a_n^{(1)}) \alpha \left(S_{\pi} ((a_1^{(2)} \ldots a_n^{(2)})^{(3)} )
(a_1^{(3)} \ldots a_{n-1}^{(3)})^{(1)}\right) \\
&& S_{\pi} \left(S_{\pi} ((a_1^{(2)} \ldots a_n^{(2)})^{(2)} )
(a_1^{(3)} \ldots a_{n-1}^{(3)})^{(2)}\right) \ot S_{\pi} ((a_1^{(2)}
\ldots a_n^{(2)})^{(1)} )\\
&&\otimes (a_1^{(3)})^{(3)} \otimes  \ldots \otimes
(a_{n-2}^{(3)})^{(3)} \beta ((a_{n-1}^{(3)})^{(3)}) \beta
(a_n^{(3)})\\
&=& \alpha (a_1^{(1)} \ldots a_n^{(1)}) \alpha \left( S (a_n^{(4)})S(a_1^{(4)}
\ldots a_{n-1}^{(4)})\, a_1^{(5)} \ldots a_{n-1}^{(5)}\right) \\
&& \sp \left( \sp (a_n^{(3)}) S(a_1^{(3)} \ldots a_{n-1}^{(3)}) \, a_{1}^{(6)} \ldots a_{n-1}^{(6)}\right) \\
&&\otimes \sp (a_{1}^{(2)} \ldots a_{n}^{(2)}) \otimes a_{1}^{(7)} \otimes \ldots \otimes a_{n-2}^{(7)}\, \beta (a_{n-1}^{(7)}) \beta (a_{n}^{(5)})\\
&=& \alpha (a_1^{(1)} \ldots a_n^{(1)}) \alpha  (S(a_n^{(4)}))\; \sp ^2
(a_n^{(3)}) \ot \sp (a_1^{(2)} \ldots a_n^{(2)})\\
&&\ot a_1^{(3)} \ot \ldots \ot a_{n-2}^{(3)} \, \beta(a_{n-1}^{(3)}) \beta({a_n^{(5)}}).
\end{eqnarray*}  Take $2 \leq j \leq n-2.$ Suppose by induction that: \begin{eqnarray*} 
t_n^j (a_1 \otimes \ldots \otimes a_n) &=& \alpha (a_1^{(1)} \ldots
a_n^{(1)}) \alpha  \left(S(a_{n-j+2}^{(4)}\ldots a_n^{(4)})\right)\; \sp ^2
(a_{n-j+2}^{(3)})\\
&& \otimes \sp ^2 (a_{n-j+3}^{(3)}) \otimes \ldots \otimes \sp ^2 (a_{n}^{(3)})  \otimes \sp(a_{1}^{(2)} \ldots a_{n}^{(2)}) \\
&&\otimes a_{1}^{(3)} \otimes \ldots \otimes a_{n-j}^{(3)} \beta
(a_{n-j+1}^{(3)}) \beta( a_{n-j+2}^{(5)} \ldots a_{n}^{(5)}). 
\end{eqnarray*} Then \begin{eqnarray*}
t_n^{j+1} (a_1 \otimes \ldots \otimes a_n) &=& \alpha  (a_1^{(1)}
\ldots a_n^{(1)}) \alpha \left(S(a_{n-j+2}^{(4)} \ldots a_n^{(4)})\right)\\
&&\alpha \left(\sp ^2 ((a_{n-j+2}^{(3)})^{(1)})  \ldots
\sp^2((a_n^{(3)})^{(1)}) \sp((a_1^{(2)} \ldots a_n^{(2)})^{(3)})
(a_1^{(3)})^{(1)} \ldots (a_{n-j}^{(3)})^{(1)}\right) \\
&& \sp \left( \sp ^2 ((a_{n-j+2}^{(3)})^{(2)})  \ldots
\sp^2((a_n^{(3)})^{(2)}) \sp((a_1^{(2)} \ldots a_n^{(2)})^{(2)})
(a_1^{(3)})^{(2)} \ldots (a_{n-j}^{(3)})^{(2)} \right) \\
&& \sp ^2((a_{n-j+2}^{(3)})^{(3)}) \ot \ldots \ot \sp
^2((a_{n}^{(3)})^{(3)}) \ot \sp((a_1^{(2)}\ldots a_n^{(2)})^{(1)}) \\
&& (a_1^{(3)})^{(3)} \ot \ldots \ot  (a_{n-j-1}^{(3)})^{(3)} \,
\beta((a_{n-j}^{(3)})^{(3)}) \beta(a_{n-j+1}^{(3)}a_{n-j+2}^{(5)}
\ldots a_{n}^{(5)})\\
&=&   \alpha (a_1^{(1)} \ldots
a_n^{(1)}) \alpha \left(S(a_{n-j+2}^{(8)} \ldots a_n^{(8)})\right)\\
&&\alpha \left(S^2
(a_{n-j+2}^{(5)} \ldots a_n^{(5)}))
 \alpha(S(a_1^{(4)} \ldots
a_n^{(4)})a_1^{(5)} \ldots a_{n-j}^{(5)}\right) \\
&& \sp \left(\sp ^2 (a_{n-j+2}^{(6)} \ldots a_n^{(6)}) \sp (a_1^{(3)}
\ldots a_n^{(3)})a_1^{(6)} \ldots a_{n-j}^{(6)}\right) \\
&& \sp ^2 (a_{n-j+2}^{(7)}) \ot \ldots \ot \sp ^2 (a_n^{(7)}) \ot \sp
(a_1^{(2)}\ldots a_n^{(2)})\\
&& \ot a_1^{(7)} \ot \ldots \ot
a_{n-j-1}^{(7)} \, \beta(a_{n-j}^{(7)}) \beta(a_{n-j+1}^{(5)})
\beta(a_{n-j+2}^{(9)} \ldots a_n^{(9)})\\%\end{eqnarray*} \begin{eqnarray*}
%t_n^{j+1} (a_1 \otimes \ldots \otimes a_n)
&=&  \alpha (a_1^{(1)} \ldots
a_n^{(1)}) \alpha \left(S(a_{n-j+2}^{(8)} \ldots a_n^{(8)})\right)\\
&& \alpha \left(S^2
(a_{n-j+2}^{(5)} \ldots a_n^{(5)}))  \alpha(S(a_1^{(4)} \ldots
a_n^{(4)})a_1^{(5)} \ldots a_{n-j}^{(5)}\right) \\
&& \sp \left(\pi S^2 (a_{n-j+2}^{(6)} \ldots a_n^{(6)}) S(a_{n-j+2}^{(3)}
\ldots a_n^{(3)})S(a_{n-j+1}^{(3)}) S(a_1^{(3)} \ldots a_{n-j}^{(3)})a_1^{(4)} \ldots
a_{n-j}^{(4)}\right)\\
&&\sp ^2 (a_{n-j+2}^{(7)}) \ot \ldots \ot \sp ^2 (a_n^{(7)}) \ot \sp
(a_1^{(2)}\ldots a_n^{(2)})\\
&& \ot a_1^{(7)} \ot \ldots \ot
a_{n-j-1}^{(7)} \, \beta(a_{n-j}^{(7)}) \beta(a_{n-j+1}^{(5)})
\beta(a_{n-j+2}^{(9)} \ldots a_n^{(9)})\\
&=&  \alpha (a_1^{(1)} \ldots
a_n^{(1)}) \alpha \left(S(a_{n-j+2}^{(4)} \ldots a_n^{(4)}))
\alpha(S(a_{n-j+1}^{(4)})\right)\\
&& \sp ^2(a_{n-j+1}^{(3)}) \ot \sp ^2
(a_{n-j+2}^{(4)}) \ot \ldots \ot \sp ^2 (a_n^{(4)}) \ot \sp
(a_1^{(2)}\ldots a_n^{(2)})\\
&& \ot a_1^{(3)} \ot \ldots \ot
a_{n-j-1}^{(3)} \, \beta(a_{n-j}^{(3)}) \beta(a_{n-j+1}^{(5)})
\beta(a_{n-j+2}^{(5)} \ldots a_n^{(5)});
\end{eqnarray*} this proves the induction. Finally, \begin{eqnarray*}
t_n^{n+1} (a_1 \otimes \ldots \ot a_n) &=& t_n^2 \left(t_n^{n-1}(a_1 \ot
\ldots \ot a_n) \right) \\
&=&t_n^2 \{ \alpha(a_1^{(1)}\ldots a_n^{(1)}) \alpha \left(S(a_3^{(4)}\ldots
a_n^{(4)}) \right) \; \sp ^2(a_3^{(3)}) \ot \sp ^2 (a_4^{(3)}) \ot \ldots \\&& \ot
\sp ^2 (a_n^{(3)}) \ot \sp (a_1^{(2)} \ldots a_n^{(2)}) \ot a_1^{(3)}
\beta(a_2^{(3)})\beta(a_3^{(5)} \ldots a_n^{(5)})\} \\
&=& \alpha(a_1^{(1)} \ldots a_n^{(1)}) \alpha \left(S(a_3^{(4)} \ldots
a_n^{(4)})\right)\\
&& \alpha \left( \sp ^2 ((a_3^{(3)})^{(1)})\ldots \sp
^2((a_n^{(3)})^{(1)}) \sp ((a_1^{(2)})^{(3)} \ldots (a_n^{(2)})^{(3)})
(a_1^{(3)})^{(1)} \right) \\
&&\alpha(S((a_1^{(3)})^{(4)})) \; \sp ^2 ((a_1^{(3)})^{(3)})\\
&& \sp \left( \sp ^2 ((a_3^{(3)})^{(2)}) \ldots \sp ^2  ((a_n^{(3)})^{(2)})
\sp((a_1^{(2)})^{(2)} \ldots (a_n^{(2)})^{(2)})(a_1^{(3)})^{(2)}
\right)\\
&&\sp ^2((a_3^{(3)})^{(3)}) \ot \ldots \ot \sp ^2((a_n^{(3)})^{(3)})
\\
&&\beta \left(\sp ((a_1^{(2)})^{(1)} \ldots (a_n^{(2)})^{(1)})  \right)
\beta((a_1^{(3)})^{(5)}) \beta(a_2^{(5)}) \beta(a_3^{(5)} \ldots
a_n^{(5)}))\end{eqnarray*} \begin{eqnarray*}t_n^{n+1} (a_1 \otimes \ldots \ot a_n)
&=&  \alpha(a_1^{(1)} \ldots a_n^{(1)}) \alpha \left(S(a_3^{(4)} \ldots
a_n^{(4)})\right)\\
&&\alpha \left(S^2(a_3^{(5)}\ldots a_n^{(5)})S(a_3^{(4)}\ldots
a_n^{(4)})S(a_2^{(4)})S(a_1^{(4)})a_1^{(5)}\right) \alpha(S(a_1^{(8)})) \;
\sp ^2 (a_1^{(7)})\\
&&\ot \sp \left(\pi S^2 (a_3^{(6)}\ldots a_n^{(6)})S(a_3^{(3)}\ldots
a_n^{(3)})S(a_2^{(3)})S(a_1^{(3)})a_1^{(6)}\right)\\
&&\ot \sp ^2 (a_3^{(7)}) \ot \ldots \ot \sp ^2 (a_n^{(7)})\;\beta \left(S(a_1^{(2)}\ldots a_n^{(2)}) \right)\beta(a_1^{(9)}a_2^{(5)}\ldots
a_n^{(5)})\\
&=&  \alpha(a_1^{(1)} \ldots a_n^{(1)}) \alpha \left(S(a_3^{(4)} \ldots
a_n^{(4)})\right) \alpha(S(a_2^{(4)}))\alpha(S(a_1^{(4)})) \; \sp ^2
(a_1^{(3)})\\
&& \ot \sp ^2 (a_2^{(3)}) \ot \ldots \ot \sp ^2 (a_n^{(3)}) \;
\beta \left(S(a_1^{(2)} \ldots a_n^{(2)})\right) \beta(a_1^{(5)}a_2^{(5)} \ldots
a_n^{(5)})\\
&=&  (\alpha * \sp * \beta )^2 (a_1) \otimes \ldots \otimes (\alpha * \sp * \beta )^2 (a_n),
\end{eqnarray*} which is equal to $ a_1 \otimes \ldots \otimes a_n$ by
condition \textsc{(cm)}, thereby  proving the result. $\square$ \end{proof}

\begin{ex} When $H$ is the trivial Hopf algebra $k,$ then $\hcm (k) $ is easy to compute, and it is equal to the classical cyclic homology $HC_* (k)$ (the datum $(\pi, \alpha, \beta)$ is necessarily equal to $(1,id_k, id_k)$).
\end{ex}

\begin{ex} Suppose $H$ is a group algebra $kG.$ Then for any element $\pi$  in the centre of $G,$ and any characters $\alpha$ and $\beta$ on $H,$ provided they take the value 1 at $\pi,$ we can consider the homology $\hcm(kG).$ We shall study examples of this situation in more detail shortly.
\end{ex}

\

Since $C_*^{(\pi, \alpha, \beta)}(H)$ is a cyclic module, we can apply  the general theory for cyclic modules, there is a long periodic  exact sequence, involving $\hcm (H)$ and the homology of the underlying simplicial module of $C_*^{(\pi, \alpha, \beta)}(H);$  view $k$ as an $H-$bimodule via: $$a. \lambda .b = \beta (a) \lambda \alpha (b), \ \forall a, b \in H \mbox{ and } \lambda \in k.$$ Let $ _{\beta}  k_{\alpha}$ denote this module. Then this last homology is in fact the Hochschild homology of the underlying algebra of $H$ with coefficients in $\kk.$ Therefore,

\begin{prop}\label{SBI} There is a long exact sequence: $$\cdots \longrightarrow H_n(H,\,  \kk)  \longrightarrow HC_n^{\pi, \alpha, \beta} (H) \longrightarrow HC_{n-2}^{\pi, \alpha, \beta} (H) \longrightarrow H_{n-1} (H,\, \kk) \longrightarrow \cdots.$$
\end{prop}

We are now going to consider the case of a group algebra.

%%%%%%%%%%%%%%%%%%%%%%%%%%%%%%%%%%%%%%%%%%%%%%%%%%%%%%%%%%%%%%%%%%%%%%%%%%%%%%%%%%%%%%%%%%%%%%%%%%

\section{Cyclic homology of a group algebra} 
\subsection{Case of trivial characters}

In this section, $H$ is a group algebra $kG,$ and the characters are both equal to the counit $\varepsilon.$ D. Burghelea, M. Karoubi and O.E. Villamayor (see~\cite{B}, \cite{KV}, \cite{L}) have established a decomposition of the classical cyclic homology of a group algebra $HC_*(kG).$ We are going to interpret this in terms of the cyclic homology of Connes and Moscovici. 

First, let us write the cyclic module maps (associated to $(\pi, \varepsilon, \varepsilon))$ for $kG:$
\begin{eqnarray*}
d_0(g_0 \otimes \ldots \otimes g_n)&=& g_1 \otimes \ldots \otimes g_n, \\
d_i(g_0 \otimes \ldots \otimes g_n)&=& g_0 \otimes \ldots \otimes g_{i-1} g_i \otimes \ldots \otimes g_n \mbox{ for $1 \leq i \leq n,$}\\
d_{n+1}(g_0 \otimes \ldots \otimes g_n)&=& g_0 \otimes \ldots \otimes g_{n-1}, \\
s_i(g_1 \otimes \ldots \otimes g_n) &=& g_1 \otimes \ldots \otimes 1 \otimes g_i \otimes \ldots \otimes g_n \mbox{ for $1 \leq i \leq n,$ and} \\
t_n(g_1 \otimes \ldots \otimes g_n) &=& \pi (g_1 \ldots g_n)^{-1} \otimes g_1 \otimes \ldots \otimes g_{n-1},
\end{eqnarray*} for any $g_0, \ldots g_n$ in $G.$

We can consider another cyclic module associated to $G$ and $\pi:$

\begin{defprop}[\cite{L} 7.4.4.] For any discrete group $G$ and any element $\pi$ in $G,$ we define a cyclic module $k\Gamma_. (G, \pi)$ as follows: as a set, $\Gamma_n (G, \pi)$ is the set of all $(g_0, \ldots , g_n) $ in $G^{n+1}$ such that $g_0 \ldots g_n$ is conjugate to $\pi,$ and the faces, degeneracies, and cyclic action are the usual ones:
\begin{eqnarray*} 
d_i(g_0, \ldots , g_n)& =& (g_0 , \ldots , g_i g_{i+1} , \ldots , g_n) \mbox{ for $0 \leq i < n,$} \\
d_n(g_0, \ldots , g_n)& =& (g_n g_0, g_1, \ldots, g_{n-1}), \\
s_i(g_0, \ldots , g_n)& =& (g_0, \ldots , g_i , 1 , g_{i+1} , \ldots , g_n), \mbox{ and} \\
t_n(g_0, \ldots , g_n)& =& (g_n, g_0, \ldots , g_{n-1}).
\end{eqnarray*} There is a canonical splitting of cyclic modules: $$C_.(kG) \cong \bigoplus_{\left<\pi \right> \in \left<G\right>} k(\Gamma_. (G, \pi)).$$
\end{defprop}

We shall now make the link with Connes and Moscovici's cyclic homology:

\begin{prop} The map \begin{eqnarray*}
(\theta _{\pi})_n :C_n^{(\pi, \varepsilon, \varepsilon)} (kG_{\pi}) & \longrightarrow & k(\Gamma_n (G, \pi)) \\
g_1 \ot \ldots \ot g_n & \mapsto & \pi (g_1 \ldots g_n)^{-1} \ot g_1 \ot \ldots \ot g_n
\end{eqnarray*} is a morphism of cyclic modules which induces an isomorphism on homology. Here $G_{\pi}$ is the centralizer of $\pi$ in $G.$ 
\end{prop}

\begin{proof} It is clear that $\theta _{\pi} $ is a simplicial map;  it also commutes with the cyclic operation, due to the fact that $\pi$ is central in $G_{\pi}.$

At the simplicial module level, there is a factorization of $\theta _{\pi}$ as 
$$C_.^{\pi, \varepsilon, \varepsilon} (kG_{\pi})=C_.(G_{\pi},k)
\stackrel{\psi}{\longrightarrow} C_.(G,k\left<\pi\right>)
\stackrel{\phi ^{-1}}{\longrightarrow} k(\Gamma_. (G, \pi)),$$ where
$\psi (g_1, \ldots , g_n) = (\pi;g_1, \ldots , g_n),$ and $\phi
(g_0,\ldots , g_n)=(g_1\ldots g_ng_0; g_1,\ldots,g_n).$ The module $k\left<\pi\right>$ is induced by the inclusion map $G_{\pi} \hookrightarrow G$ from the trivial $G_{\pi}-$module $k.$ Therefore Shapiro's lemma implies that $\psi$ is a quasi-isomorphism. Since $\phi$ is also a quasi-isomorphism (see~\cite{L} 7.4.2.), so is $\theta _{\pi}.$
Finally, the long periodic exact sequences for both cyclic modules, together with the five lemma, give an isomorphism of the cyclic homologies.  $\square$ \end{proof}

\begin{rmk} This proof is very close to the proof of \cite{L} 7.4.5.
\end{rmk}

These two propositions combined yield the following result:

\begin{theo} For any discrete group $G,$  there is a graded isomorphism:
$$HC_*(kG) \cong \bigoplus_{\left<\pi \right> \in \left<G\right>} HC_*^{\pi, \varepsilon, \varepsilon} (kG_{\pi}).$$
\end{theo}

\begin{rmk} \label{interp} Using the results in \cite{L} (7.4.11 to 7.4.13), there are various interpretations of $\hcme (kG_{\pi})$ which we can give:

(i) If $G$ is a torsion free group, then $$HC_n^{\pi, \varepsilon, \varepsilon} (kG_{\pi})= \begin{cases} H_n(G_{\pi}/\{\pi\},k) & \mbox{ if $\pi \neq 1$,} \\ H_n(G) \oplus H_{n-2}(G) \oplus \ldots  & \mbox{ if $\pi = 1$}.
\end{cases}$$

Note that a similar result for Lie algebras was obtained by Connes and Moscovici in \cite{CM} and \cite{CM3}, in the cohomological framework (for trivial character and grouplike).

(ii) If $G$ is abelian, then $$HC_*^{\pi, \varepsilon, \varepsilon}(kG) = \begin{cases} H_*(G/\{\pi\},k) & \mbox{ if $\pi $ is of infinite order in $G,$} \\ HC_*(k) \ot H_*(G/\{\pi\},k) & \mbox{ if $\pi $ is of finite order in $G$} \end{cases}$$ as a graded module.

(iii) If $k$ contains $\mathbb{Q},$ then $$HC_*^{\pi, \varepsilon, \varepsilon} (kG_{\pi})= \begin{cases} H_*(G_{\pi}/\{\pi\},k) & \mbox{ if $\pi $ is of infinite order in $G,$} \\ HC_*(k) \ot H_*(G_{\pi}/\{\pi\},k) & \mbox{ if $\pi $ is of finite order in $G$} \end{cases}$$ as a graded module.

In all these expressions, $\{\pi\}$ is the cyclic subgroup of $G$ generated by $\pi.$
\end{rmk}

\begin{rmk} We shall see further on, on some examples (quiver algebras), that these decomposition formulae cannot be extended to general Hopf algebras (other than group algebras).
\end{rmk}

\begin{rmk} These interpretations also enable us to compute explicitly $\hcme (kG)$ when $G$ is a cyclic group, and therefore the classical $HC_*(kG)$ also:
\end{rmk}

\begin{prop} \label{cyclgp} Let $G$ be a finite cyclic group, and let $\pi$ be an element in $G.$ Let $m_{\pi}$ be the index of $\pi$ in $G.$ Let $k$ be a ring, and view it as a trivial $k \left(\mathbb{Z} / m_{\pi}\mathbb{Z}\right) -$module. Then:$$HC_n^{\pi,\varepsilon, \varepsilon} (kG) = \begin{cases} k \oplus Ann(m_{\pi})^{n/2} & \mbox{ if $n$ is even} \\ (k/m_{\pi}k)^{n+1/2} & \mbox{ if $n$ is odd}, \end{cases}$$ where $Ann(m_{\pi})= \{ \lambda \in k / \; m_{\pi} \lambda =0 \}.$
\end{prop}

\begin{proof}  Let $\sigma _{\pi}$ be a generator for $\mathbb{Z} / m_{\pi}\mathbb{Z} ,$ and let $N_{\pi}=1+ \sigma _{\pi} + \sigma _{\pi} ^2 + \ldots + \sigma _{\pi} ^{m_{\pi}-1}$ be its norm. The homology of $\mathbb{Z} / m_{\pi}\mathbb{Z} $ with coefficients in $k$ is given as follows (see~\cite{W}): $$H_n(\mathbb{Z} / m_{\pi}\mathbb{Z} ;k)= \begin{cases} k/(\sigma _{\pi}-1)k = k & \mbox{ if $n=0,$} \\ k^G/N_{\pi} k=k/m_{\pi} k & \mbox{ if $n$ is odd,} \\ \{\lambda \in k / N_{\pi} \lambda =0\}/(\sigma_{\pi}-1)k=Ann(m_{\pi}) &  \mbox{ if $n$ is even $>0$}.\end{cases}$$  Remark \ref{interp} (ii) yields the result. $\square$ \end{proof}

\begin{cor}\label{cyclgpcor} With the same notations, we can express the classical cyclic homology of $kG:$
$$HC_n(kG)= \begin{cases} k^{\#G} \oplus \left( \bigoplus_{\pi \in G} Ann(m_{\pi}) ^{n/2} \right) & \mbox{ if $n$ is even} \\ \bigoplus_{\pi \in G}(k/m_{\pi}k)^{n+1/2} & \mbox{ if $n$ is odd}.  \end{cases}$$
\end{cor}

\begin{rmk} This agrees with the result in~\cite{BACH2}.
\end{rmk}

\begin{rmk} We can consider some special cases of these results. For instance, if $m_{\pi}$ is not a divisor of zero in $k,$ then $HC_n^{\pi, \varepsilon,\varepsilon} (kG)$ is equal to k when $n$ is even, and to $(k/m_{\pi}k)^{n+1/2}$ when $n$ is odd, so that $HC_n(kG)$ is equal to $k^{\#G}$ when $n$ is even, and  to $ \bigoplus_{\pi \in G}( k/m_{\pi}k)^{n+1/2}$ when $n$ is odd; if moreover the order of $G$ is prime, this gives Theorem 1 in~\cite{CGV}.
\end{rmk}

\begin{rmk} Note that the results of Proposition \ref{cyclgp} and
Corollary \ref{cyclgpcor} remain true for non-cyclic groups $G,$ as
long as $G/\{\pi\}$ is cyclic and $\pi $ is in the centre of $G$.
\end{rmk}

%%%%%%%%%%%%%%%%%%%%%%%%%%%%%%%%%%%%%%%%%%%%%%%%%%%%%%%%%%%%%%%%%%%%%%%%%%

\subsection{Case of non-trivial characters when $G$ is a cyclic group} 

In this paragraph, we assume that $k$ is a characteristic zero field. Set $H=kG$ with $G=\mathbb{Z}/m\mathbb{Z},$ and let $\alpha$ and $\beta$ be characters on $H.$ In the first place, we shall compute the Hochschild homology $H_*(H, \, _{\beta}k_{\alpha}).$ For this, we shall use a simplified projective resolution defined in \cite{CGV}:

$$\xymatrix{\ldots \ar[r] &  H^{\ot 2} \ar[rr]^-{\sum_{i=0}^{m-1} g^i\ot g^{m-i}} && H^{\ot 2} \ar[rr]^-{-1\ot 1+g\ot g^{m-1}} && \\
\ldots \ar[r] &  H^{\ot 2} \ar[rr]^-{\sum_{i=0}^{m-1} g^i\ot g^{m-i}} && H^{\ot 2} \ar[rr]^-{-1\ot 1+g\ot g^{m-1}} &&  H^{\ot 2}  \ar[r]^-{\mu} & H \ar[r] &0}$$ where $\mu $ is the multiplication in $H,$ $g$ is a fixed generator in $\mathbb{Z}/m\mathbb{Z},$ and the maps are multiplication by the terms above the arrows.

Suppose first that $\alpha \neq \beta.$ Tensoring the above resolution by $\kk$ over $H^e$ yields the following complex:

$$\xymatrix{\ldots \ar[rr]^{\zeta \rho ^{-1}-1} && \kk \ar[r]^{0} &
\kk \ar[rr]^{\zeta \rho ^{-1}-1}& & \ldots \ar[r]^{0} &  \kk
\ar[rr]^{\zeta \rho ^{-1}-1}&& \kk \ar[r] & 0}$$ in which $\zeta =
\alpha(g)$ and $\rho = \beta (g)$ are $m^{th}$ roots of unity in $k$. The
homology of this complex is zero:  $H_*(H, \, \kk)=0.$

The SBI exact sequence of Proposition \ref{SBI} yields the cyclic
homology of $H$ in this case. 

We shall now study the case $\alpha = \beta$ (we will not use the Hochschild homology in this case).

Let $\pi = g^s$ be an element in $G$ with $\zeta ^s =1$ (ie. $\alpha(\pi) =1).$ Consider the map 
\begin{eqnarray*}
\chi_n^{\pi, \zeta}: \, C_n^{(\pi, \alpha, \alpha)} (kG) & \longrightarrow &  C_n^{(\pi, \varepsilon, \varepsilon)} (kG) \\
g^{i_1} \ot \ldots \ot g^{i_n} & \mapsto & \zeta^{i_1 + \cdots + i_n} g^{i_1} \ot \ldots \ot g^{i_n}.
\end{eqnarray*} Then $\chi^{\pi,\zeta}_.$ is an isomorphism of cyclic modules (the inverse is $\chi^{\pi,\zeta^{-1}}_.).$

Therefore, $\hcm(kG)$ is isomorphic to $\hcme(kG).$ Using Proposition \ref{cyclgp}, we finally have:

\begin{prop} The cyclic homology of the Hopf algebra $H=k\left(\mathbb{Z}/m\mathbb{Z}\right)$ is given as follows: $$\begin{array}{l} 
\hcm(H)=0 \mbox{ if $\alpha \neq \beta,$} \\
HC_n^{\pi, \alpha, \alpha}(H) = \begin{cases}  k \oplus Ann(m_{\pi})^{n/2} & \mbox{ if $n$ is even} \\ (k/m_{\pi}k)^{n+1/2} & \mbox{ if $n$ is odd}. \end{cases}
\end{array}$$
\end{prop}

\section{Cyclic homologies of some truncated quiver algebras}

\subsection{Hochschild homologies  of truncated quiver algebras}

Let $\Delta$ be a finite quiver (that is a finite oriented graph), $k$
a commutative ring, and let $k \Delta$ be the algebra of paths in
$\Delta$ ($k\Delta$ is the free $k-$module with basis the set of paths
in $\Delta,$ and the multiplication is obtained via the concatenation
of paths). For all $p \in \mathbb{N},$ let $ \Delta _p$  denote the set
of paths of length $p$ in $\Delta,$ and let  $\frak{m}$ denote the
ideal of $k \Delta$ generated by $\Delta_1.$ Consider an admissible
ideal $I$ in $k\Delta$ (that is, an ideal $I$ such that there exists
an integer $n$ with $\frak{m}^n \subset I \subset \frak{m}^2$). In~\cite{AG}, Anick and Green introduced a new quiver $\Gamma$ which gave them a minimal projective $k\Delta /I-$resolution of the algebra $k\Delta_0 ,$ graded by the lengths of paths.
 
When $I$ is equal to $\frak{m} ^n,$ for an integer $n,$
the set of vertices of this quiver $\Gamma$  is $\Delta_0 \cup \Delta _1 \cup \ldots \cup \Delta _{n-1}, $ and the edges are given as follows:
$$\begin{array}{cl}
a \leftarrow e & \mbox{if } a \in \Delta_1, e \in \Delta_0 \mbox{ and the terminus of $a$ is $e,$} \\
\gamma \leftarrow a & \mbox{if } a \in \Delta_1, \gamma \in \Delta_{n-1} \mbox{ and the terminus of $\gamma$ is the origin of $a,$ } \\
a \leftarrow \gamma  & \mbox{if } a \in \Delta_1, \gamma \in \Delta_{n-1} \mbox{ and the terminus of $a$ is the origin of $\gamma.$}
\end{array}$$
If $\Gamma ^{(i)}$ denotes the set of paths of length $i$ in $\Gamma,$ then $\Gamma ^{(2c)}$ can be identified with $\Delta_{nc}$ and $\Gamma ^{(2c+1)}$ can be identified with $\Delta_{nc+1}.$ 

\subsubsection{First case: $n\geq2$}

Using Anick and Green's resolution, Sköldberg, in~\cite{S}, constructed a resolution of the algebra $A := k\Delta / \frak{m}^n$ for $n \geq 2$ as follows:

\begin{theo}[\cite{S} Theorem 1] \label{res}There is a projective $A-$bimodule resolution of $A,$ which is graded by the length of paths, as follows:
$$\mathbf{P}_A: \cdots \stackrel{\tiny d_{i+1}}{\longrightarrow} P_i \stackrel{\tiny d_{i}}{\longrightarrow} \cdots  \stackrel{\tiny d_{2}}{\longrightarrow} P_1\stackrel{\tiny d_{1}}{\longrightarrow} P_0  \stackrel{\tiny d_{0}}{\longrightarrow} A \longrightarrow 0, $$ where $$P_i = A \otimes _{k\Delta _0} k\Gamma ^{(i)} \otimes _{k\Delta _0} A,$$  the differentials are defined by 
$$\begin{array}{l}
d_{2c+1}( u  \otimes a_1 \cdots  a_{cn+1} \otimes v)   = u  a_1 \otimes a_2 \cdots a_{cn+1} \otimes v \ - \ u \otimes a_1 \cdots a_{cn} \otimes a_{cn+1}  v ,\\
\\
d_{2c}( u \otimes  a_1 \cdots a_{cn}  \otimes v)\\
\\  
\ \  = \sum_{j=0}^{n-1} u a_1 \cdots a_j \otimes a_{j+1} \cdots a_{(c-1)n+j+1} \otimes a_{(c-1)n+j+2} \cdots a_{cn} v \  \mbox{ if $c > 0,$} \\
\end{array}$$
and the augmentation $d_0 : A \otimes _{k\Delta _0} A \cong  A \otimes _{k\Delta _0} k\Gamma ^{(0)} \otimes _{k\Delta _o} A \longrightarrow A$ is defined by $$d_0 ( u \otimes v ) = u v .$$

\end{theo}

Note that the differentials preserve the gradation, so that the Hochschild homology spaces also are graded. Let  $HH_{p,q} (A)$ denote the $q^{th}$ graded part of the space $HH_p (A).$

By means of the resolution in Theorem \ref{res},  Sköldberg computes the Hochschild homology of $A$ with coefficients in $A;$ to state this result,  we shall need some notation: let $\mathcal{C}$ denote the set of cycles in the quiver $\Delta,$ and for any cycle $\gamma$ in $\mathcal{C},$ let $L(\gamma)$ denote its length. There is a natural action of the cyclic group $\left< t_{\gamma} \right>$ of order $L(\gamma)$ on $\gamma;$ let $\overline{\gamma}$ denote the orbit of $\gamma$ under this action, and let  $\overline{\mathcal{C}}$ denote the set of orbits of cycles.

\begin{theo} [\cite{S}]\label{S}  Set $q=cn+e,$ for $0 \leq e
\leq n-1 \ (n \geq 2).$ Then the Hochschild homology space $HH_{p,q}
(A)$ is given by:
$$ \begin{cases} 
k^{a_q} & \mbox{if  $1 \leq e \leq n-1$  and $ 2c \leq p \leq 2c+1,$} \\
\bigoplus _{r|q} (k^{(n\wedge r)-1 } \oplus Ker (. \frac{n}{n\wedge r} : k \rightarrow k))^{b_r} & \mbox{if  $e=0,$  and $ 0<2c=p,$} \\
\bigoplus _{r|q} (k^{(n\wedge r)-1 } \oplus Coker (. \frac{n}{n\wedge r} : k \rightarrow k))^{b_r} & \mbox{if $ e=0,$  and $ 0<2c-1=p,$} \\
k^{\# \frak{C}_0} & \mbox{if $ p=q=0,$  and} \\
0 & \mbox{otherwise,}
\end{cases}  $$ where $a_q$ is the number of cycles of length $q$ in
$\overline{\mathcal{C}},$   $b_r$ is the number of cycles of length
$r$ in $\overline{\mathcal{C}}$ which are not powers of smaller
cycles, and $n \wedge r$ is the greatest common divisor of $n$ and
$r.$
\iffalse Set $q=cn+e,$ for $0 \leq e \leq n-1 \ (n \geq 2).$ Then:
$$HH_{p,q} (A) = \begin{cases} 
k^{a_q} & if \ 1 \leq e \leq n-1 \ and \ 2c \leq p \leq 2c+1, \\
\bigoplus _{r|q} (k^{(n\wedge r)-1 } \oplus Ker (. \frac{n}{n\wedge r} : k \rightarrow k))^{b_r} & if \ e=0, \ and \ 0<2c=p, \\
\bigoplus _{r|q} (k^{(n\wedge r)-1 } \oplus Coker (. \frac{n}{n\wedge r} : k \rightarrow k))^{b_r} & if \ e=0, \ and \ 0<2c-1=p, \\
k^{\# \Delta_0} & if \ p=q=0, \ and \\
0 & otherwise,
\end{cases}  $$ where $a_q$ is the number of cycles of length $q$ in $\overline{\mathcal{C}},$  and $b_r$ is the number of cycles of length $r$ in $\overline{\mathcal{C}}$ which are not powers of smaller cycles, and $n \wedge r$ is the greatest common divisor of $n$ and $r.$\fi
\end{theo}

\begin{ex} Suppose $\Delta = \Delta^{(n)}$ is the $n-$crown, that is the quiver with $n$ vertices $e_0, \ldots , e_{n-1},$ and  $n$ edges $a_0, \ldots , a_{n-1},$ each edge $a_i$ going from the vertex $e_i$ to the vertex $e_{i+1}$ for $0 \leq i \leq n-2,$ and the edge $a_{n-1}$ going from $e_{n-1}$ to $e_0,$ as follows:$$\epsfysize4cm\epsfbox{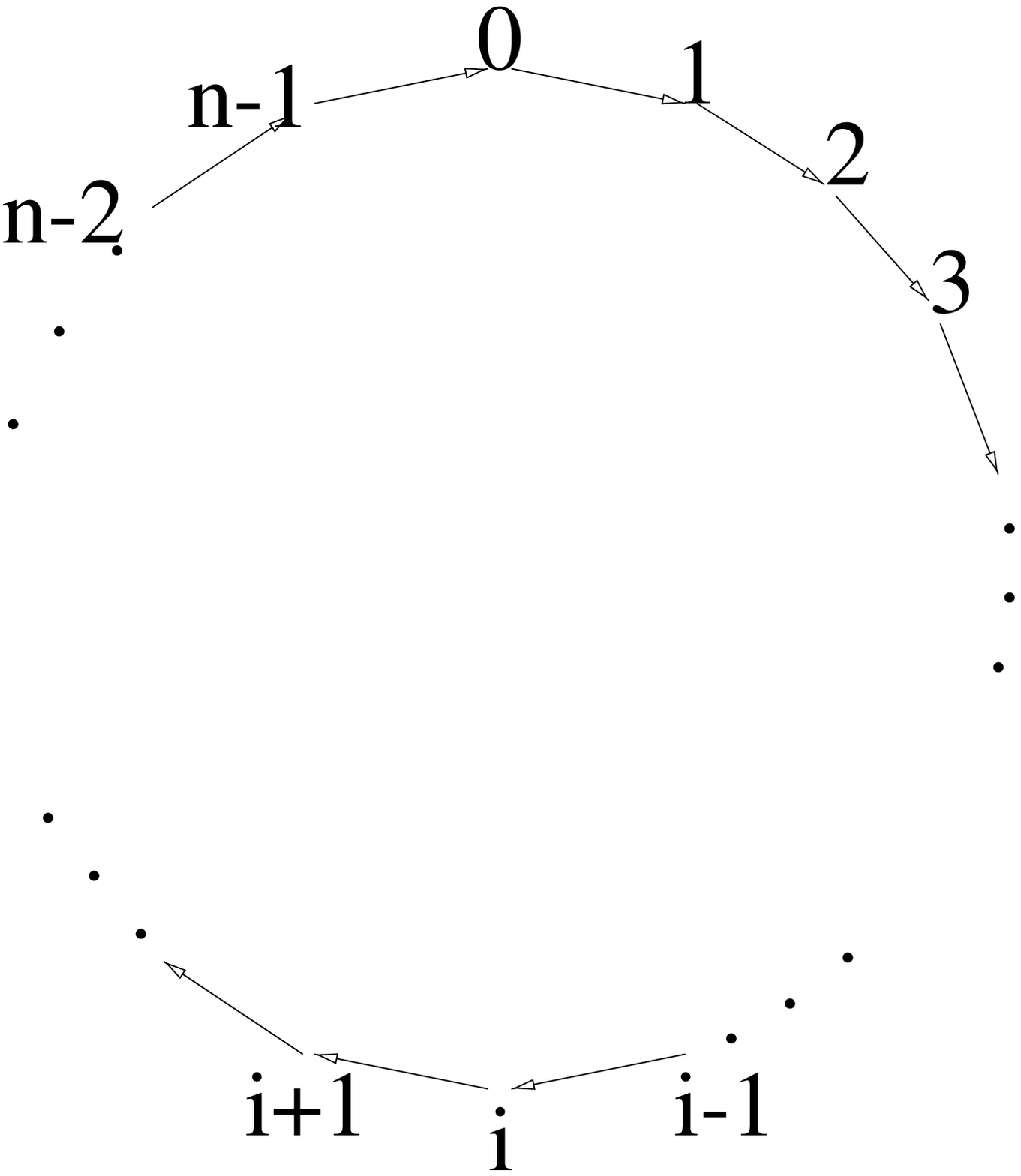}$$
 Then the Hochschild homology of the Taft algebra $\Lambda _n :=k\Delta / \frak{m}^n$ is given by:
\begin{eqnarray*}
HH_{p,cn} (\Lambda _n )& =&k^{n-1} \mbox{ if $p=2c$ or $p=2c-1$}\\
HH_{0,0} (\Lambda _n ) &=& k^n \\
HH_{p,q} (\Lambda _n ) &=& \ 0 \mbox{ in all other cases}.
\end{eqnarray*} \end{ex}

\begin{rmk} The cyclic cohomology of these truncated quiver algebras has been computed in \cite{BLM} and \cite{Li}.
\end{rmk}

\

The resolution of Theorem~(\ref{res}) also enables us to calculate other Hochschild homologies, useful for the computation of $\hcm (A).$ Let $\alpha $ and $\beta$ be characters on $A;$  we are going to compute the Hochschild homology $H_* (A, \, _{\beta}  k_{\alpha}),$ for $n \geq 2.$  This is equal to $Tor_*^{A^e} (A,\,\kk),$ that is to the homology of the complex $ _{\beta}  k_{\alpha} \otimes _{A-A} {\bf P} _A.$ This complex is isomorphic to the complex $ _{\beta}  k_{\alpha} \otimes_{k \Delta_0^e} \Gamma ^{(\bullet)}$ with the following differentials: $$\begin{array}{rcl}
 _{\beta}  k_{\alpha} \otimes_{k \Delta_0^e} \Gamma ^{(2c)} & \longrightarrow & _{\beta}  k_{\alpha} \otimes_{k \Delta_0^e} \Gamma ^{(2c-1)} \\
1 \otimes b_1 \ldots b_{nc} & \mapsto & \sum_{j=0}^{n-1} \alpha (b_{(c-1)n+j+2}  \ldots b_{cn}) \beta (b_1 \ldots b_j) \otimes b_{j+1} \ldots b_{(c-1)n +j+1} \\
\\
 _{\beta}  k_{\alpha} \otimes_{k \Delta_0^e} \Gamma ^{(2c+1)} & \longrightarrow & _{\beta}  k_{\alpha} \otimes_{k \Delta_0^e} \Gamma ^{(2c)} \\
1 \otimes b_1 \ldots b_{nc+1} & \mapsto &  \alpha(b_1)\otimes b_2 \ldots b_{nc+1} - \beta (b_{nc+1}) \otimes b_1 \ldots b_{nc}.
\end{array}$$

Note that the differentials do not preserve the gradation in this case.

Now the list of characters on $A$ is the following: for each vertex $e_i,$ there is one character $\alpha _i$ which equals 1 on $e_i,$ and 0 on every other path in $\Delta$ (the $e_i$ are orthogonal idempotents, and the edges either have different extremities, or have a power which is zero). In particular, the characters vanish on all the paths of length greater than 1, so the differentials vanish. Therefore, the homology is:

\begin{eqnarray*} 
H_{2c} (A,\, _{\beta}  k_{\alpha})& =&  _{\beta}  k_{\alpha} \otimes _{k \Delta_0^e} k \Delta _{nc} \\
H_{2c+1} (A,\, _{\beta}  k_{\alpha})& =&  _{\beta}  k_{\alpha} \otimes _{k \Delta_0^e} k \Delta _{nc+1}. 
\end{eqnarray*}

\begin{ex}  For the Taft algebras, the results are:

$$\begin{cases} 
H_{2c} (\Lambda _n , \, _{\beta}  k_{\alpha} ) =k \\
H_{2c+1} (\Lambda _n , \, _{\beta}  k_{\alpha} ) =0  
\end{cases}\ if \ \alpha = \beta,$$
$$\begin{cases} 
H_{2c} (\Lambda _n , \, _{\beta}  k_{\alpha} ) =0 \\
H_{2c+1} (\Lambda _n , \, _{\beta}  k_{\alpha} ) =k  
\end{cases}\ if \  \alpha =  \alpha _{i+1} \  and \ \beta = \alpha _i ,$$ and in all other cases the homology is 0 in all degrees.\end{ex}
 
\subsubsection{Second case: $n=0$ or $n=1$}

We shall now look at the cases $n=0$ and $n=1.$ 
The case $n=1$ is quite simple: $k\Delta/\m$ is equal to $k \Delta _0 \cong \times _{s \in \Delta_0} ks,$ so that $$HH_p(k\Delta /\m)=\bigoplus _{s \in \Delta_0} HH_p(ks) = \begin{cases} \bigoplus _{s \in \Delta_0} ks & \mbox{ if $ p=0,$} \\ 0 & \mbox{ if $p>0.$} \end{cases}$$

Finally, for the case $n=0,$ we can  state:

\begin{prop} The Hochschild homology of $k\Delta$ is given by: $$\begin{cases} HH_0 (k\Delta)=k\overline{\mathcal{C}}, \\ HH_1 (k\Delta)= \{\sum_{i=0}^{L(\gamma) -1}t_{\gamma}^i (\gamma)/\ \gamma \in \overline{\mathcal{C}}, L(\gamma) \geq 1\} \\  HH_p (k\Delta)=0 \mbox{ if $p \geq 2$}. \end{cases}$$
\end{prop}

\begin{proof} We shall use the following resolution  (see for instance~\cite{C2}):
\begin{lemma}[\cite{C2} Theorem 2.5] \label{resC} There is a $k\Delta-$bimodule projective resolution of $k\Delta$ given by $$\ldots  0\longrightarrow k\Delta \ot _{k\Delta_0} k\Delta_1 \ot_{k\Delta_0} k\Delta \longrightarrow  k\Delta \ot _{k\Delta_0} k\Delta \longrightarrow  k\Delta  \longrightarrow 0. $$ \end{lemma}
Tensoring by $k\Delta$ over $k\Delta^e$ yields the following complex: $$\begin{array}{ccccc}
\ldots \longrightarrow 0 \longrightarrow &  k\Delta \ot _{k\Delta_0^e} k\Delta_1& \stackrel{\delta}{\longrightarrow} &  k\mathcal{C} & \longrightarrow 0. \\
&\nu \ot a & \mapsto & \nu a - a \nu & \end{array}$$ 

The space $HH_0(k\Delta)$ is generated by the cycles in $\mathcal{C},$ subjected to the relations given by the image of $\delta.$ Since $\delta (\nu \ot a)=\nu a - t_{\nu a }(\nu a ),$ the relations identify two cycles in the same orbit, and  $HH_0(k\Delta)=k\overline{\mathcal{C}}.$

The complex is $\overline{\mathcal{C}}-$graded; therefore, to find $HH_1(k\Delta)=\ker \delta,$ it is sufficient to consider elements of type $x=\sum_{i=0}^{L(\gamma)-1} \lambda_i t_{\gamma}^i (\gamma)$, in which $\nu \ot a$ is identified with $\nu a,$ and the $ \lambda_i$ belong to $k.$ We have $\delta (x) =0$ iff $ \lambda_0 = \lambda_1= \ldots = \lambda_{L(\gamma)-1},$ and the result follows. $\square$ \end{proof}

\begin{rmk} We may there again compute the homology $H_*(A, \, \kk),$ when $n$ is equal to 0 or 1.

For $n=0,$ we shall once more use the resolution of lemma~\ref{resC}. Tensoring it by $\kk$ over $k\Delta ^e,$ we obtain the following complex: $$\ldots 0 \longrightarrow \, \kk \ot_{k\Delta _0^e} k\Delta _1 \longrightarrow \, \kk \ot_{k\Delta _0^e} k\Delta _0 \longrightarrow 0.$$  All the maps in this complex are zero (the characters vanish on the edges). Given any character $\chi,$ let $e_{\chi}$ denote the unique vertex such that $\chi (e_{\chi})=1.$ Then $H_0(k\Delta,\, \kk) = \, \kk \ot_{k\Delta _0^e} k\Delta _0 \cong e_{\alpha} k\Delta _1e_{\beta},$ and $H_1(k\Delta) =  \kk \ot_{k\Delta _0^e} k\Delta _0  $ is equal to 0 if $\alpha \neq \beta,$ and to $ke_{\alpha}$ if $\alpha = \beta.$

For $n=1,$ we have $H_*( ks,\, \kk) =H_*(k,k)$ if $\alpha = \beta$ or if both characters differ from $\alpha_s;$ in the other cases, $H_*( ks,\, \kk) =0. $ 

Hence, $H_0(k\Delta/\m,\, \kk)=k^{\#\Delta_0}$ if $\alpha = \beta, \ k^{\#\Delta_0 -2}$ if $\alpha \neq \beta, $  and $H_p(k\Delta/\m,\,\kk)=0$ for all $p>0.$
\end{rmk}

%%%%%%%%%%%%%%%%%%%%%%%%%%%%%%%%%%%%%%%%%%%%%%%%%%%%%%%%%%%%%%%%%%%%%%%%%%%%%%%%%%%%%%%%%%%%%%%%%%%%%%%%%%%

\subsection{Cyclic homology of graded algebras}

In this paragraph, $k$ is a commutative ring which contains $\mathbb{Q}.$ When $A$ is a graded $k-$algebra, Connes' SBI exact sequence splits in the following way:

\begin{theo} [\cite{L} Theorem 4.1.13]\label{SBIgr} Let $A$ be a unital graded algebra over $k$ containing $\mathbb{Q}.$ Define $\overline{HH}_p (A) = HH_p(A) / HH_p (A_0)$ and $\overline{HC}_p (A) = HC_p(A) / HC_p (A_0).$ Connes' exact sequence for $\overline{HC}$ reduces to the short exact sequences: $$0 \rightarrow \overline{HC}_{n-1}(A) \rightarrow \overline{HH}_n (A)\rightarrow \overline{HC}_n(A) \rightarrow 0.$$  
\end{theo}

This will enable us to compute the classical cyclic homology of truncated quiver algebras. Let us first consider the cases $n=0$ and $n=1.$ Combining the results for Hochschild homology and Theorem~\ref{SBIgr} yields the following:

\begin{prop}\label{cycl01} The cyclic homology of $k\Delta$ and of $k\Delta/\m$ are given by: \begin{eqnarray*} 
HC_{2c}(k\Delta /\m)&=& \oplus_{s\in \Delta_0} ks\\
HC_{2c+1}(k\Delta /\m)&=& 0 \\
&\mbox{and}\\
HC_{0}(k\Delta)&=& k\overline{\mathcal{C}} \\
HC_{2c}(k\Delta )&=& k^{\#\Delta_0} \\
HC_{2c+1}(k\Delta )&=& 0,
\end{eqnarray*} for all $c \in \mathbb{N}.$
\end{prop}

The case $n\geq 2$ involves the same methods:

\begin{prop} Suppose $n \geq 2.$ Then: 
\begin{eqnarray*}
\mathrm{dim}_k \,HC_{2c} (k\Delta / \m ^n) &=& \#\Delta _0 + \sum_{e=1}^{n-1} a_{cn+e} - \sum_{\tiny \begin{array}{c} r|(c+1)n \\% r | \hspace{-3pt} / n 
n \notin r \mathbb{N} \end{array}} (r \wedge n -1)b_r \\ 
\mathrm{dim}_k \,HC_{2c+1} (k\Delta / \m ^n) &=&  \sum_{r|n} ( r -1) b_r.
\end{eqnarray*}
\end{prop}

\begin{proof} In the first place, $A_0$ is equal to $k\Delta_0,$ so that we know the homologies of  $A_0$ (see Proposition \ref{cycl01}). Next, we have $HC_0(A) = HH_0(A) = k^{\#\Delta_0 + \sum_{e=1}^{n-1} a_e}.$ Then, using Theorem \ref{SBIgr}, we get the following formula: $$dim_k HC_{2c}(A) +dim_k HC_{2c+1}(A)=\#\Delta_0 + \sum_{r|(c+1)n}(r \wedge n -1)b_r + \sum_{e=1}^{n-1} a_{cn+e}.$$ In particular, $dim_k HC_{1}(A) = \sum_{r|n}(r-1)b_r.$

An induction on $c$ yields the result. $\square$ \end{proof}

\begin{cor} When $\Delta$ is the $n-$crown, then the results are:
$$\begin{cases} HC_{2c} (\Lambda _n ) = k^n, \\
HC_{2c+1} (\Lambda _n ) = k^{n-1} \mbox{ for $c \in \mathbb{N}.$}
\end{cases}$$
\end{cor}

\begin{ex} Let $A$ be the quotient algebra $k[X]/(X^n).$ It has a presentation by quiver and relations (the quiver has one vertex and one loop), and its cyclic homology is $k^n$ in even degree, and vanishes in odd degree.
\end{ex}

\begin{rmk} This agrees with the general results  given in~\cite{BACH} and~\cite{BACH2} (in which $k$ may be a field of positive characteristic, or in fact any commutative ring). 
\end{rmk}
 
%%%%%%%%%%%%%%%%%%%%%%%%%%%%%%%%%%%%%%%%%%%%%%%%%%%%%%%%%%%%%%%%%%%%%%%%%%

\subsection{Connes and Moscovici homology of some truncated algebras}

Here, $k$ is a  commutative ring which contains a primitive $n-$th root of unity $q.$ The Taft algebra $\lan$ is then a Hopf  algebra (see~\cite{C1}), with the following structure maps:
\begin{eqnarray*} \varepsilon (e_i)=\delta _{i,0}, &&  \varepsilon (a_i)=0, \\
\Delta (e_i)= \sum_{j+k=i} e_j \ot e_k, && \Delta (a_i)= \sum_{j+k=i}( e_j \ot a_k+q^k a_j \ot e_k), \\
S(e_i)=e_{-i}, && S(a_i) = -q^{i+1} a_{-i-1}, 
\end{eqnarray*} where $\delta$ is the Kronecker symbol.
We can therefore consider the homology $\hcm (\Lambda _n).$

First of all, we need to find out which characters and grouplikes $(\pi, \alpha, \beta)$ satisfy the necessary conditions \textsc{(cm)}. There are $n$ grouplike elements in $\Lambda _n ,$ given as follows: $$\pi _i = \sum_{l=0}^{n-1} q^{il}e_l.$$ A triple $(\pi _i,\alpha _u, \alpha _v)$ satisfies \textsc{(cm)} iff $ui \equiv 0 (\mathrm{mod} \ n), \ vi \equiv 0 (\mathrm{mod} \ n),$ and $v-u+1+i \equiv 0 (\mathrm{mod} \ n).$ For instance, when $u=v,$ this means that $i$ is necessarily equal to -1, and $u$ and $v$ are equal to 0 (that is $\alpha_u=\varepsilon=\alpha_v$); when $v=u-1,$ it means that $i$ is equal to 0, that is, $\pi _i$ is equal to 1. We shall not need to know the details of the other possibilities.

Using the long periodic exact sequence for $\hcm$ (Proposition \ref{SBI}), we obtain the following results: 

\begin{prop} 
\begin{eqnarray*}
HC_p^{\pi_{n-1}, \varepsilon, \varepsilon} (\Lambda _n) =& \begin{cases} k^{p/2 +1} & \mbox{ if $p$ is even,} \\ 0 &\mbox{ if $p$ is odd,} \end{cases} \\
HC_p^{1, \alpha_u, \alpha_{u-1}} (\Lambda _n) =& \begin{cases} 0 & \mbox{ if $p$ is even,} \\ k^{p+1/2 }  &\mbox{ if $p$ is odd,} \end{cases} & \mbox{ for all $u \in \{0, \ldots, n-1 \},$}\\
\hcm (\Lambda _n) =& 0 & \mbox{ in every other case.}
\end{eqnarray*}
\end{prop}

\begin{rmk} Note that none of these homologies are direct factors in the classical $HC_* (\Lambda _n),$ which seems to preclude any possibility of decomposing $HC_* (\Lambda _n)$ as a sum of Connes and Moscovici homologies.
\end{rmk}

%\begin{thebibliography}{99}

%\bibitem[AG]{AG} \textsc{Anick, D. and Green, E. L.,} On the homology of quotients of path algebras, \textit{Comm. Alg.} \textbf{15} (1987), pp 309-341. \\

%\vspace*{-15pt}

%\bibitem[C]{C} \textsc{Cibils, C.,} A quiver quantum group, \textit{Comm. Math. Phys.,} \textbf{157} (1993), no. 3, pp 459-477. \\

%\vspace*{-15pt}

%\bibitem[L]{L} \textsc{Loday, J-L.,} Cyclic homology, Appendix E by María O. Ronco, \textit{Springer-Verlag, Berlin,} (1992), xviii+454 pp. ISBN:3-540-53339-7 \\

%\vspace*{-15pt}

%\bibitem[S]{S} \textsc{Skoldberg, E.,} The Hochschild homology of truncated and quadratic monomial algebras, \textit{J. London Math. Soc. (2)} \textbf{59} (1999), no.1, pp 76-86. \\

%\vspace*{-15pt}

%\end{thebibliography}

%\end{document}

\end{document}